\documentclass{article}
\usepackage{float}
\usepackage[caption = false]{subfig}
\usepackage[final]{graphicx}

\usepackage{amsthm, amsmath, amssymb, amsfonts}
\usepackage{mathtools}
\usepackage{tikz}
\usetikzlibrary{calc}
\usetikzlibrary{patterns.meta}

\usepackage{bbm}
\usepackage{url}
\usepackage{cite}

\usepackage[utf8]{inputenc}
\usepackage{eurosym}
\usepackage{framed}
\usepackage{enumerate}
\usepackage[linesnumbered, boxed, lined, ruled,vlined]{algorithm2e}

\usepackage{todonotes}

\newcommand{\RL}{\text{SRC}}
\newcommand{\RA}{\text{SRE}}
\newcommand{\RM}{\text{aFRR}}
\newcommand{\DA}{\text{DA}}

\newcommand{\TC}{\text{PT}}

\title{Optimizing the Marketing of Flexibility for a Virtual Battery in Day-Ahead and Balancing Markets: A Rolling Horizon Case Study}
\author{E. Finhold \and C. Gärtner \and R. Grindel \and T. Heller \and N. Leithäuser \and E. Röger \and F. Schirra}

\date{Fraunhofer Institute for Industrial Mathematics ITWM\\  67663 Kaiserslautern\\  Germany}

\begin{document}

\maketitle

\begin{abstract}
Industrial electricity consumers with flexible demand can profit by adjusting their load to short-term prices and by providing balancing services to the grid. Markets which support this kind of short-term position adjustment are the day-ahead market and balancing markets. We propose a formulation for a combined optimization model that computes an optimal distribution of flexibility between the balancing and day-ahead markets. The optimal solution also includes the specific bids for the day-ahead and balancing markets. Besides the expected profits of each market and their individual bidding languages, our model also takes their different roles in a continuous marketing of flexibility into account. To prevent overrating short-term profits we introduce a variable penalty term that adds a cost to unfavorable load schedules.
We evaluate the optimization model in a rolling horizon case study based on the setting of a virtual battery at TRIMET SE, which is derived from a flexible aluminum electrolysis process. For such a battery we compute a daily optimal split of flexibility and trading decisions based on data in the period 04/2021 - 03/2022. We show that the optimal split is more profitable than using only one market or a fixed split between the markets.
\end{abstract}


\section{Introduction}\label{sec: intro}
Regenerative electricity production in Germany has increased considerably over the last two decades. With this transformation away from schedulable fossil fuels to mostly volatile regenerative power sources, electricity production becomes more uncertain. As a consequence, the consumption of electricity is also bound to change: it has to become more flexible in order to cope with the uncertainty in renewable energy production. Electricity consumers should consume more if there is an oversupply of (renewable) energy and less if energy is scarce. 
This flexibility in energy consumption is referred to as a \emph{virtual battery}, as instead of physically storing energy, a similar effect is achieved by shifting the demand. 

Since electricity prices are not constant, neither within a day nor over the course of several days, a flexible consumer can arbitrage trade the difference between times with higher prices against those with lower prices. Usually, these price differences reflect the fluctuations in the availability  of electricity. Thus, consumers who schedule consumption based on prices can simultaneously reduce their own costs and support the shift to renewable energy sources. The general term for this approach of minimizing electricity costs is \emph{demand side management} (DSM). For a general introduction to DSM, we refer to the surveys by Zhang and Grossmann (cf. \cite{zhang2016enterprise}), Siano (cf. \cite{siano2014demand}), and Gelazanskas and Gamage (cf. \cite{gelazanskas2014demand}). Applications of DSM can be found in several kind of industries, such as steel production and aluminium electrolysis (cf. \cite{schafer2019model}). 

The day-ahead market is of particular importance for marketing flexibility or closing short-term positions (cf. \cite{epexspot}). Trading decisions for day-ahead markets can either be specified by an amount for a given hour, e.g. as in \cite{leithauser2022optimal, kraft2022stochastic, boomsma2014cmnordic}, or by exploiting a more complex bidding language, e.g. by using piece-wise linear bidding functions (cf. \cite{faria2011day, ayon2017aggregators}). Next to the day-ahead market, balancing markets allow pre-qualified providers another profit opportunity (cf. \cite{entsoe}). Optimal bidding on balancing markets requires an adaptation to specific market designs, such as price limits or minimum bid quantities. These design parameters have been changed several times since November 2020 for the German balancing markets. A bidding model and evaluation on the period from January 2014 to December 2016 can be found in \cite{ocker2018bidding}. For a more recent evaluation of a bidding optimization model on the period from April 2021 to November 2021, see \cite{gartner2022optimal}. In this article, we focus on the period from April 2021 to March 2022. Note that the price limit of 9.999,99~\euro/MWh was removed during this period, and thus, we consider a period in which the market design changes. Information concerning balancing markets is published at \cite{regelleistung}.

Due to this complexity, it is non-trivial for market participants to decide how to optimally market flexibility on a particular energy market. Obviously, it is even harder to decide how to spread flexibility across multiple markets, especially in advance. An overview over several studies on profits obtained in different energy markets with different market and price assumptions can be found in \cite{klaeboe2013optimal}. Note that combining the bidding over more than one market typically results in higher profits.

Optimization models for bidding across multiple markets are often formulated as multi-stage optimization problems (cf. \cite{boomsma2014cmnordic, bohlayer2018demand, bohlayer2020energy, lohndorf2022value, ottesen2018multi}). Thereby, different combinations of markets are studied, e.g. day-ahead market and intraday market (e.g. \cite{ayon2017aggregators, lohndorf2022value}), day-ahead market and balancing markets (e.g. \cite{boomsma2014cmnordic, schafer2019model}) and all three of these markets (e.g. \cite{kraft2022stochastic, Nolzen2022}). The intraday market is often modeled as if it were a single auction and, thus, structurally identical to the day-ahead market. Here, we neglect the intraday market and focus on the combined bidding on the day-ahead market and the balancing markets, since these markets only need a single trading decision. In contrast to the multi-stage approach, the trading on multiple markets can also be modeled as simultaneous decision as in \cite{Nolzen2022, kraft2022stochastic}. In this article, we are interested in a split of flexibility between markets, and thus, follow the simultaneous approach for our optimization model. 

When bidding on multiple energy markets, a key aspect is to avoid myopic profit realizations that are less profitable in the long run. For this, we introduce a rolling horizon framework which is based on price forecasts and additional constraints in the optimization problem. In \cite{leithauser2022optimal}, different lookahead horizons were tested for day-ahead trading, whereas Corinaldesi et al. (cf. \cite{CORINALDESI2020100392}) evaluated their trading decisions on the day-ahead and intraday market in a rolling horizon setting. In this article, we present a mathematical optimization model for a cross-market situation considering the day-ahead and the balancing market in a rolling-horizon setting. The model returns the locally optimal split of flexibility capacities for the different markets, together with the concrete bids to place on the markets. Yet, the model design ensures that the solutions are also profitable long-term. Note that in this approach, we execute only one optimization run per day where we fix all decision variables and hence do not follow a multi-stage approach that uses the realizations of one auction for the bidding in a later auction (as e.g. in \cite{ottesen2018multi}). One particular advantage of this strategy is that it allows to split the flexibility budgets among separate financial management departments specialised in different markets.

Our setting is inspired by the situation at TRIMET Aluminium SE, where the production processes have been redesigned such that the aluminium ovens can act as a large virtual battery with positive and negative capacity. In general, our approach can be transferred to common batteries. Note however that we implicitly assume a very large storage capacity with respect to the hourly flexibility and that we do not take (un)charging cycles into account. Trading decisions on single or multiple markets have already been discussed and described in various contexts, for example for flexibility aggregators (e.g. \cite{ayon2017aggregators, ottesen2018multi}), for energy storage owners (e.g. \cite{lohndorf2022value}) and for energy intensive chemical plants (e.g. \cite{XENOS2016418}). In \cite{schafer2019model}, besides a model for combined bidding strategies for the primary balancing market and the day-ahead market, the authors provide a case study for energy intense processes in the aluminium production at TRIMET Aluminium SE. While the general idea is similar to our approach in principal, the authors focus on a different balancing market (primary instead of secondary balancing market) than we do. Also, the market regulations have changed significantly in the last years. Pure day-ahead trading decisions for this setting at TRIMET Aluminium SE have also been discussed earlier by the authors of this work for single days in \cite{finhold2023bicriteria} and in a rolling horizon context in \cite{leithauser2022optimal}.

Another purpose of our study is to assess the profit opportunities of demand side management in general. We consider the problem of using flexibility profitably on the day-ahead and balancing markets. We distinctly put our focus on optimizing the profits of a single provider of flexibility. The flexibility we market stems from a flexibilized industrial process. We therefore consider not only forecasts that predict prices or profit opportunities in electricity markets, but also technical constraints such as battery size and production efficiency losses. Our key results are: 
\begin{itemize}
    \item We adapt our optimization to the market design in the 2021-2022 period and capture the intricacies of the balancing market by modelling both the power and the reserve part of the auction. 
    \item We evaluate this optimization strategy by computing cumulative profits in the rolling horizon period from April 2021 to March 2022. This rolling horizon setting plays a major role in real-world applications. Therefore, we put a special focus on it in the numerical evaluations.
    \item We show that our methods perform well in a quickly changing market and illustrate the rising value of flexibility in periods with high prices and volatility. Besides that, we see that the day-ahead market and the balancing markets fulfill different roles in the rolling horizon context. 
    
\end{itemize}
The remainder of the article is structured as follows. In Section~\ref{sec: setting}, we discuss for which trading decisions our model can provide suitable decision support and explain the \emph{rolling horizon setting} for our case. In Section~\ref{sec: markets}, we give a brief description of the considered markets, namely, the balancing markets and the day-ahead market. Here, we also provide an explanation on how to compute the required input data for our model. In Section~\ref{sec: cm problem}, we introduce our decision support model for \emph{cross-market trading} where we consider trading strategies for the balancing markets and the day-ahead market. Numerical results are then presented in Section~\ref{sec: numerical case study}. Finally, we conclude in Section~\ref{sec: outlook} with a short outlook on future research. 

\section{Decision framework}\label{sec: setting}
In the following, we describe our setting. First, we give an overview of the electricity markets considered for trading and when the trading decisions are made in the respective markets. We then briefly describe the \emph{rolling horizon}, which is an essential component for practical usage. 

\subsection{Trading decisions}\label{subsec: trading decision}

German electricity trading takes place on several independent markets. On the one hand, there is the futures and derivatives market, which can be used to build up positions at an early stage. Here, trading decisions take place up to six years before delivery of the electricity.  Consequently, they are based on predictions that cannot include any information about the situation close to delivery. Therefore, in addition to this long-term market, there are markets for short-term adjustment of existing positions, when more recent data is available. These include: The day-ahead market (DA), where trading decisions can be made one day before delivery, the intraday market (ID), which allows trading from the previous day up to five minutes before delivery, and the balancing energy market (secondary reserve or \textit{automatic Frequency Restoration Reserves} (\RM)), which is used if a grid imbalance occurs despite the adjustment of positions beforehand. It consists of a reserve capacity and a reserve energy market that are denoted by \RL~market and \RA~market, respectively. Obviously, short-term markets play a crucial role in promoting flexibility. Therefore, we concentrate on trading decisions on these markets in the following.

\begin{figure}
    \centering
    \includegraphics[width=12cm]{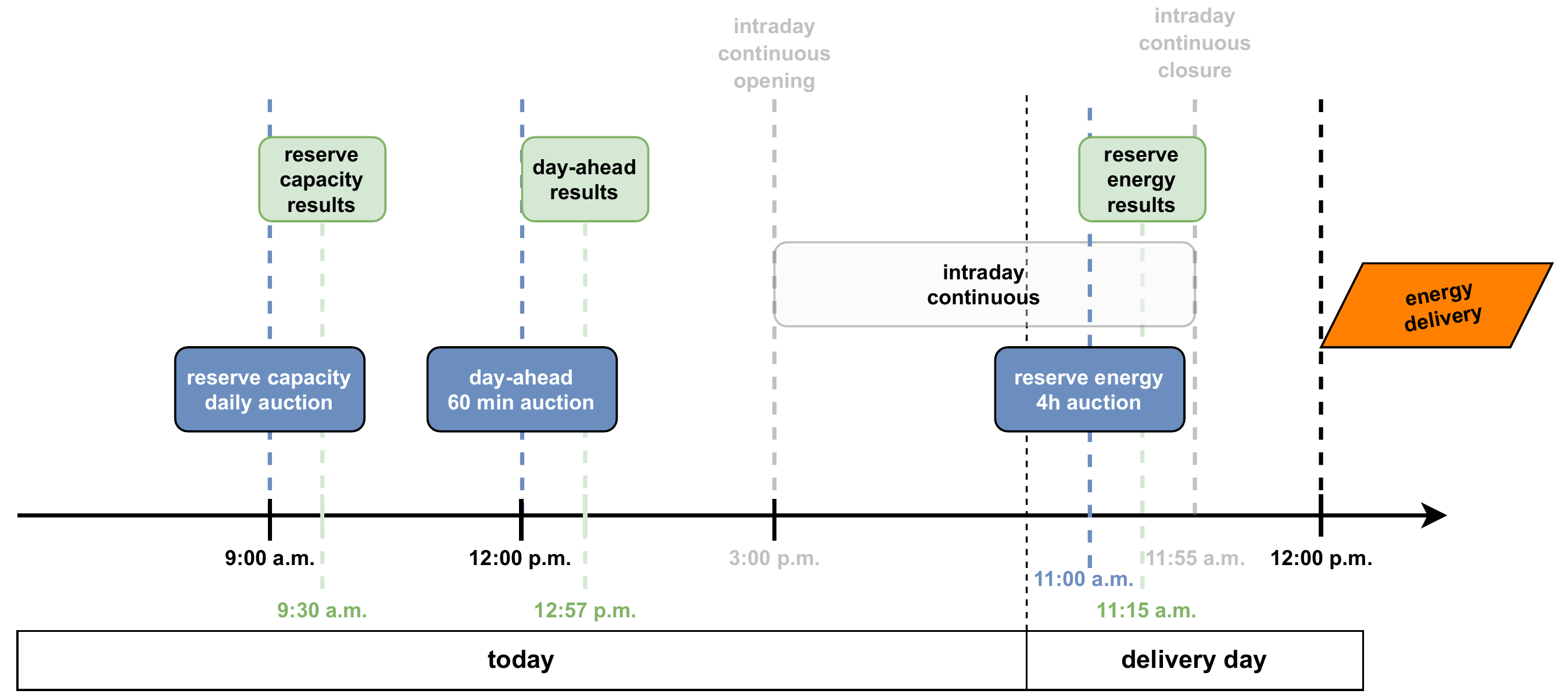}
    \caption{Auctions and trading time horizons for the German short-term electricity markets.}
    \label{fig: tm trading decision schema}
\end{figure}

The time horizons of decisions on the short-term markets is shown in Figure \ref{fig: tm trading decision schema}, exemplary for a product with delivery starting at 12~p.m. Because the data used in this article is from the years 2021 and 2022, all times are expressed in a way they fit into the market description of that period. The first decision to be made today, i.e. the day before delivery, is the one regarding volume and pricing of offers on the \RL~market; they must be submitted before 9~a.m. Results are published half an hour later. Then, the next auction that takes place is the day-ahead auction at 12~p.m. Here, results are published 57 minutes later. These auction times are independent of the product that is considered. In contrast to this, order books for \RA~bids close one hour before delivery and their results are published 15 minutes later. Consequently, this could take place today or on the delivery day itself, depending on the considered product.
Note that on the \RA~market, only six products are sold, each covering four hours of the day. The corresponding auction is always conducted one hour before the delivery period starts. Finally, the continuous intraday market starts at 3~p.m. and remains open until 5 minutes prior to delivery of each product. In the figure, this is depicted in grey, because we do not consider this market further in this article. This is due to the fact that trading on the intraday market requires either a person or an algorithm to monitor prices and to react to them. Not every company that would have access to the intraday market in theory can use it in practice due to these requirements. Therefore, we concentrate on a setting in which flexibility is split between the day-ahead and \RM~markets, that is relevant for smaller companies. As a result of both considered markets, traders must decide the level of flexibility they want to provide across all three markets before the end of the \RL~bid collection period.

\subsection{The rolling horizon}\label{subsec: rolling horizon}
Computing an optimal bidding strategy for a single day might lead to unstable plans and shortsighted decisions. Instead, each day we compute a bidding strategy for a longer look-ahead period. However, the strategy will only be executed for the first day of the period and recomputed on the following day when more recent forecasts become available. Like this. we get a setting in which the bidding strategy adapts each day to new information about the price forecasts.  

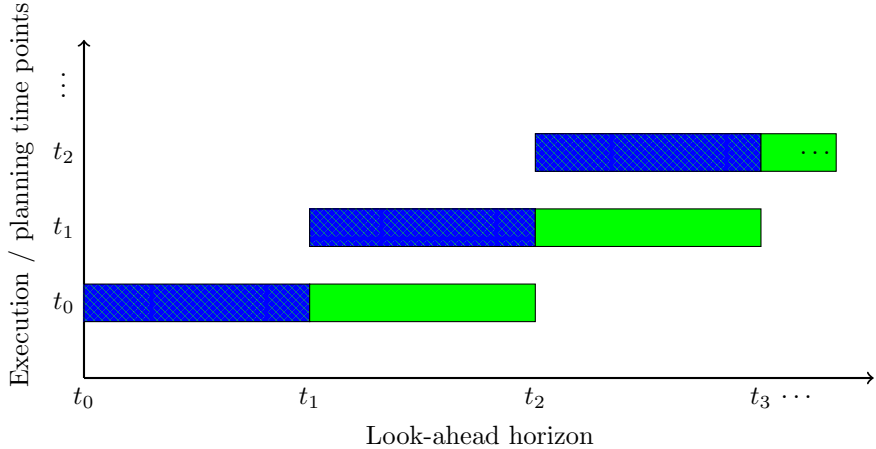
\begin{figure}
\centering
\begin{tikzpicture}

    \draw[->, thick] (0,0) node (start) [below] {$t_0$} -> (10.5,0) node (end) {};
    \draw (3,0) node [below] {$t_1$};
    \draw (6,0) node [below] {$t_2$};
    \draw (3,-0.75) node[label={[label distance=0.5cm,text depth=-1ex]right:Look-ahead horizon}] {};
    
    \draw (9,0) node [below] {$t_3$};
    \draw (9.5,-0.1) node [below] {$\dots$};
    \draw[->, thick] (0,0) -> (0, 4.5);
    \draw (0, 1) node [left] {$t_0$}; 
    \draw[pattern={Lines[angle=45, distance={3pt/sqrt(2)}, line width={5pt}]}, pattern color=blue,
    preaction={%
    pattern={Lines[angle=45, distance={3pt/sqrt(2)}, line width={5pt}, xshift=1.5pt]}, 
    pattern color=green}] (0, 0.75) rectangle ++(3, 0.5);
    \filldraw [fill=green, draw=black] (3, .75) rectangle ++(3, 0.5);
    \draw (0, 2) node [left] {$t_1$};
    \draw[pattern={Lines[angle=45, distance={3pt/sqrt(2)}, line width={5pt}]}, pattern color=blue,
    preaction={%
    pattern={Lines[angle=45, distance={3pt/sqrt(2)}, line width={5pt}, xshift=1.5pt]}, 
    pattern color=green}] (3, 1.75) rectangle ++(3, 0.5);
    \filldraw [fill=green, draw=black] (6, 1.75) rectangle ++(3, 0.5);
    \draw (0, 3) node [left] {$t_2$};

    \draw[pattern={Lines[angle=45, distance={3pt/sqrt(2)}, line width={5pt}]}, pattern color=blue,
    preaction={%
    pattern={Lines[angle=45, distance={3pt/sqrt(2)}, line width={5pt}, xshift=1.5pt]}, 
    pattern color=green}] (6, 2.75) rectangle ++(3, 0.5);
    \filldraw [fill=green, draw=black] (9, 2.75) rectangle ++(1, 0.5);
    
    \draw (9.75, 3) node {$\dots$};
    \draw (-0.1, 4) node [left] {$\vdots$};
    \draw (-1, -0.75) node[label={[label distance=0.5cm,text depth=-1ex,rotate=90]right:Execution / planning time points}] {};
\end{tikzpicture}
\caption{Schematic representation of the rolling horizon setting with a look-ahead horizon of two time points. On the x-axis, we denote the time points for which a plan is executed or computed. On the y-axis, we denote the time points where a part of a computed plan is executed and a new plan is computed. The solid green part depicts the planning horizon at the given point in time whereas the blue striped part depicts the execution horizon at the given point in time.} \label{fig: rolling horizon}
\end{figure}

To describe the rolling horizon setting in more detail, we denote the points in time at which decisions have to be made by $t_0 ,t_1, \ldots$. Let $T$ be the number of days in our planning horizon. At $t_0$, an optimal schedule for the entire planning horizon is computed and executed for the following day. At the next time~$t_1$, the procedure is repeated. Regardless of whether new price forecasts are available, a new optimal schedule is computed for the entire planning horizon. Thus, at each point in time, a schedule covering $T$ days is computed and the first day of the schedule is executed. A schematic representation of the rolling horizon can be seen in Figure~\ref{fig: rolling horizon}. 

In this article, we propose a hybrid approach where we consider the balancing market only for the next day, while the day-ahead market is considered over a longer period. For a more detailed discussion of the length of the horizon considered for bidding in the day-ahead market, we refer to \cite{leithauser2022optimal}. Note that since in general price forecasts change over time, the resulting trading schedule in a rolling horizon framework is usually not optimal, although the bidding problem is solved optimally every day.

\section{Markets and their data}\label{sec: markets}
In the following, we describe the considered markets on which we want to market the flexibility. For each market, we describe the required input data. 

\subsection{The balancing markets}\label{subsec: balancing markets}
In Germany, the electricity grid has to maintain a electric frequency of 50Hz in order to keep the electricity system stable. Transmission system operators are responsible to maintain these 50Hz by using balancing energy when the frequency moves away from its set value. This balancing energy is sold on three markets, namely the primary reserve market, the secondary reserve market and the tertiary reserve market. The main difference between the markets is the time span within which the respective products can be activated: Primary reserve products need to be available within 5 seconds after being requested, secondary reserve products take over after 30 seconds, and finally tertiary reserve products cover the range of 12.5 minutes up to 60 minutes after the demand request. To be able to handle frequency fluctuations on both directions, all three markets offer the possibility to buy and sell products with positive or negative balancing energy.

In the following paragraphs, we present the secondary reserve markets in more detail. On both markets - \RL~and \RA - pre-qualified reserve energy providers are allowed to trade. For this pre-qualification, certain criteria need to be fulfilled. These contain elements as admission to EPEX, forecasting electricity input every day, or having a minimum capacity of 5~MW and ensure that all basic necessities for offering \RM~energy are fulfilled. After that, a second pre-qualification process has to be passed where data about the power plant is transmitted to the transmission system operators. Then, a market participant is allowed to offer electricity products on the \RM~market.

All \RM~products can provide either positive or negative reserve energy and consist of 4h-blocks, during which the promised electricity needs to be available. Furthermore, as mentioned above, it needs to be available within 30 seconds after the demand request. After that, activating all balancing energy may take up to 4.5 minutes; then, all demanded energy needs to be fully powered and has to be available for 15 minutes in total. We denote the set of 4h block products on the balancing markets as $\mathcal{I}$. It consists of intervals~$I_k=[a_k,b_k)$ spanning the time period of block $I_k$, whereas the corresponding index set is given by $\mathcal{K}=\{1,\dots, 6\}$ and the products are indexed in their natural order, i.e. $I_1 = [0, 4), ..., I_6 = [20, 0)$. In the following, we reference product $I_k$ by $k$ if the meaning is obvious through the context.

Regarding regulations in the \RM~market, the minimum bidding amount was 5~MW until July 2018; since then, this has been reduced to 1~MW. Nonetheless, this change comes with the restriction that, when 1~MW to 4~MW are offered on the market, then only one offer per product~$I_k \in \mathcal{I}$ is allowed. All offers are collected daily and consist of an electricity amount, whether positive or negative reserve energy is offered, and two prices: One price belongs to the \RL~market and represents the remuneration for holding electricity available. The other price belongs to the \RA~market and is paid for actually activated electricity. If a market participant's \RL~market price was too high and the offer is not awarded, its \RA~market price is nonetheless placed in the merit order curve of the \RA~market. Consequently, also those market participants are allowed to offer electricity on the \RA~market that did not make it in the \RL~market.

\subsubsection{The reserve capacity market}\label{subsubsec: reserve capacity}
On the \RL~market, remuneration for offered electricity is determined by a pay-as-bid method. All offers are collected in a merit order based on their \RL~prices. The merit order is used to determine which offers are awarded. When an offer is awarded, its price is paid to the market participant for holding electricity available, independent of whether it will be used or not.

To depict this procedure in our model, we define \emph{bidding scenarios} on the auction of the \RL~market. Each bidding scenario~$i$ is described by a price level~$p_i^{\RL}$ in \euro/MW and a volume~$v_i^\RL$ in MW. In order to keep the model linear, we restrict the possible bidding prices to the discrete set of the mentioned price levels. Consequently, we consider $N^{\RL}$ bidding scenarios that are in total ordered by increasing price levels, i.e. $p_i^\RL \leq p_{i+1}^\RL$ for all $i=1, \dots, N^{\RL}-1$. 

Corresponding to these bidding scenarios, we denote by $q_i^\RL$ the probability that exactly the first $i$ bids are awarded. Since the price levels are ordered increasingly, price level~$p_i^{\RL}$ is the highest accepted price level. Consequently, $q_0^{\RL}$ denotes the probability of no accepted price level. Note that this approach can also be found in \cite{kraft2022stochastic}. 
The calibration of all introduced parameters are based on historical data. In order to obtain the probabilities~$q_{k, i}^{\RL}$ that a bid $p_i^{\RL}$ for product~$k$ is accepted on the \RL~market, we use the time series of daily marginal prices per product~$k$. In the following, we describe how to obtain said probabilities for a fixed day. Let the marginal accepted price for product~$k$ at day~$d$ be denoted as $\pi_{d,k}^{\RL}$. To get an estimate of the short-term distribution of the marginal prices and to strip ourselves from the assumption of stationarity of the time series, the calibration period is restricted to the time interval of 30 days before the considered day. Consequently, we obtain an empirical distribution function for each day~$d$ as
\begin{align}
    \hat{F_d}(p) = \frac{1}{30} \sum_{j = d-31}^{d-1}\mathbbm{1}_{\pi_{j,k}^{\RL} < \, p} \;,
\end{align}
where the summation goes over the amount of times during the 30 day interval that the marginal accepted price was lower than a threshold~$p$. Consequently, the $q_i^\RL$ for a day $d$ can be computed by 
\begin{align*}
    q_{k, 0}^{\RL} &= \hat{F_d}(p_1^\RL),  \\
    q_{k, 1}^{\RL} &= \hat{F_d}(p_2^\RL) - \hat{F_d}(p_1^\RL),\\
    q_{k, 2}^{\RL} &= \hat{F_d}(p_3^\RL) - \hat{F_d}(p_2^\RL),  \\
    & \vdots\\
    q_{k, N}^{\RL} &= 1 - \hat{F_d}(p_{N^\RL}).
\end{align*}

\begin{figure}[ht]
    \centering
    \includegraphics[width=12cm]{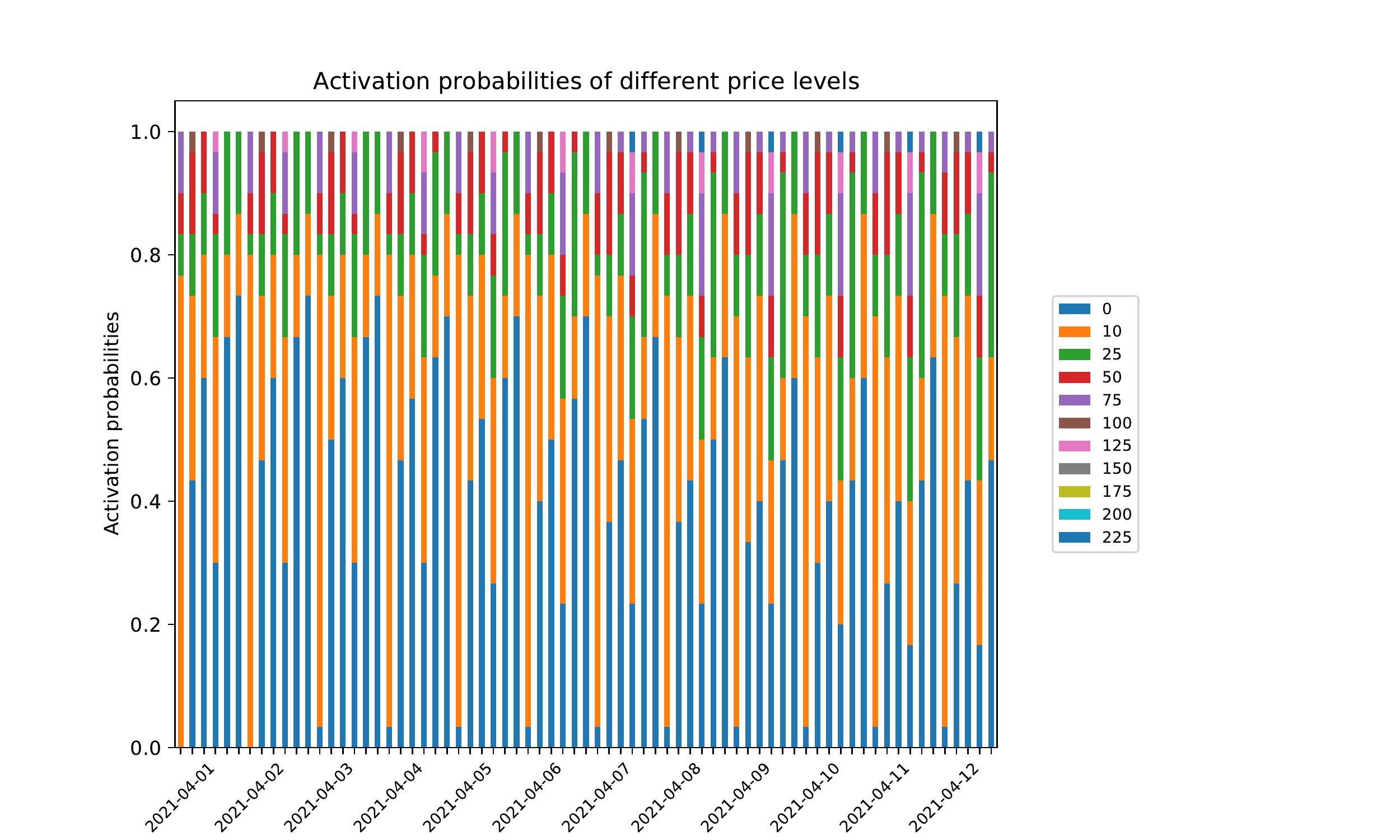}
    \caption{Probabilities for different price levels over time, where each price level is associated with one color. }
    \label{fig:probabilitiespricelevels}
\end{figure}

For illustration purposes, we show in Figure~\ref{fig:probabilitiespricelevels} the increasingly ordered cumulative probabilities for different price levels. For example, the first bar in the bar chart describes the first product of 01/04/2021. There, the price~0 has a probability of zero of being the highest accepted bid. In the shown period, the price level of~225 is the highest price level with a probability strictly greater than zero. Over the course of a day, a pattern can be observed. Although the calculated probability of an acceptance of price 0 is higher for later products on the day, the highest accepted price is also higher. We will use these price probabilities later on as input for the daily optimization runs. 

\subsubsection{The reserve energy market}\label{subsubsec: reserve energy}

Following the same market design as the \RL~market, the \RA~market is also based on the pay-as-bid system in combination with a merit order. The only difference is that the offers' \RA~market prices are used instead of \RL~market prices. Furthermore, next to the \RM~bids, so-called \textit{free bids} can be entered or adapted for the \RA~market until gate closure. Free bids are placed by market participants that did not bid on the \RL~market or whose bids were not accepted, but who still want to participate in the \RA~market.

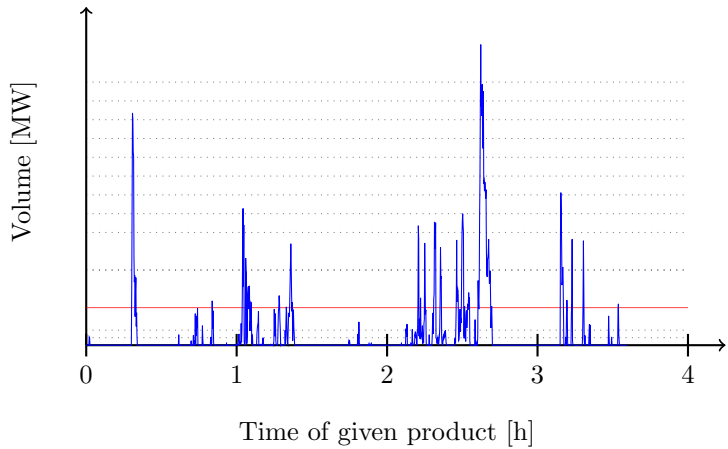
\begin{figure}
    \centering
    \begin{tikzpicture}
    % axis
        \draw[->, thick] (0,0) node (start) [below] {} -> (8.5,0) node (end) {};
        \draw (-1, 0.75) node[label={[label distance=0.5cm,text depth=-1ex,rotate=90]right:Volume [MW]}] {};
            \draw (4, -0.25) node[label={[label distance=0.5cm,text depth=-1ex]below:Time of given product [h]}] {};
        \draw (0,-0.15) node [below] {$0$};
        \draw[-, thick] (0,0.15) -> (0,-0.15);
        \draw (2,-0.15) node [below] {$1$};
        \draw[-, thick] (2,0.15) -> (2,-0.15);
        \draw (4,-0.15) node [below] {$2$};
        \draw[-, thick] (4,0.15) -> (4,-0.15);
        \draw (6,-0.15) node [below] {$3$};
        \draw[-, thick] (6,0.15) -> (6,-0.15);
        \draw (8,-0.15) node [below] {$4$};
        \draw[-, thick] (8,0.15) -> (8,-0.15);
        \draw[->, thick] (0,0) -> (0, 4.5);

    % bids
        \draw[dotted, opacity=0.5] (0,1) -> (8,1);
        \draw[dotted, opacity=0.5] (0,0.1) -> (8,0.1);
        \draw[dotted, opacity=0.5] (0,0.2) -> (8,0.2);
        \draw[dotted, opacity=0.5] (0,1) -> (8,1);
        \draw[dotted, opacity=0.5] (0,1.5) -> (8,1.5);
        \draw[-, red, opacity=0.7] (0,0.5) -> (8,0.5);
        \draw[dotted, opacity=0.5] (0,2) -> (8,2);
        \draw[dotted, opacity=0.5] (0,1.75) -> (8,1.75);
        \draw[dotted, opacity=0.5] (0,2.25) -> (8,2.25);
        \draw[dotted, opacity=0.5] (0,2.5) -> (8,2.5);
        \draw[dotted, opacity=0.5] (0,2.75) -> (8,2.75);
        \draw[dotted, opacity=0.5] (0,3) -> (8,3);
        \draw[dotted, opacity=0.5] (0,3.25) -> (8,3.25);
        \draw[dotted, opacity=0.5] (0,3.5) -> (8,3.5);
        
    % activation function   
    \draw[blue, thin] plot[] file {Activation_Calls.txt};
    \end{tikzpicture}
    \caption{Schematic representation of the computation of the activation duration for the positive \RA~market. On the x-axis, the 4h of one product are shown. On the y-axis, we denote the offered and demand volumes on the \RA~market. The horizontal lines represent the cumulative volume of bids on the \RA~market, which are ordered by their corresponding prices. The red line is the bid of the considered price level, whereas the blue line is a schematic representation of the activation function.}
    \label{fig: activation duration}
\end{figure}

On the \RL~market, the price for holding energy available is determined for every market participant through the procedure described above. In contrast to this, the demanded quantity of an offer depends on the actual amount of \RM~energy that is needed. This can change from second to second. Consequently, the amount that is paid for is the aggregated time during which electricity was actually activated times the price. This \emph{activation amount} in MW is denoted by $S_k(t)$ for product~$I_k$ and for every second~$t$ of the product time interval. Furthermore, we define $\Psi_k(p)$ as the sum of offered capacity in MW at the \RA~market for product~$k$, aggregated over all offers with an ask price of at most $p$. This contains the amount of capacity available for product~$I_k$ up to a threshold price~$p$. With this, we are able to define the \emph{\RL~activation duration} for a given product~$I_k$ and a price level~$p$ as
\begin{align}\label{SRL_act_duration}
    L_k(p) = \int_{I_k}\mathbbm{1}_{S_k(t)\geq \Psi_k(p)}dt. 
\end{align}
In Figure~\ref{fig: activation duration}, an example for the activation duration and its computation is given. Whenever the blue lines lies above the red line in that figure, the considered bid is activated. Note that the considered price levels for the \RA~market are not necessarily the same as the price levels for the \RL~market.

Now, the \emph{expected activation duration}~$\alpha_{k, i}$ of \RA~activation for a given price offer~$p_i$ and a product~$k$ is defined by
\begin{align}\label{SRL_exp_act_duration}
    \alpha_{k, i} \coloneqq \mathbb{E}\left[L_k(p_i)\right]. 
\end{align}

In order to calculate the expected activation duration~$\alpha_{k, i}$, we first compute the \RA~activation duration~$L_k(p_i)$ as in \eqref{SRL_act_duration} for the given price levels~$p_1, \dots, p_{N^{\RA}}$ for every day in the data set and use then the seven-day average as the expected value. 

\subsection{The day-ahead market}\label{subsec: da market}
Here we focus on the German day-ahead market on the European Power Exchange (EPEX SPOT SE). The day-ahead market is still the most important electricity market for spot products, as can be seen from the volumes traded on this market (cf. \cite{epexspot}). This market is conducted through an auction where each hour of the following day is traded as a single hourly product. All participants, both buyers and sellers, pay a single market clearing price with an additional spread. The auction design allows for complex bids, but the strategy we focus on is simply to bid high/low enough to be accepted at all reasonable market clearing prices. For the detailed rules of the market we refer to \cite{epexspot}. 

\begin{figure}
    \centering
    \includegraphics[width=12cm]{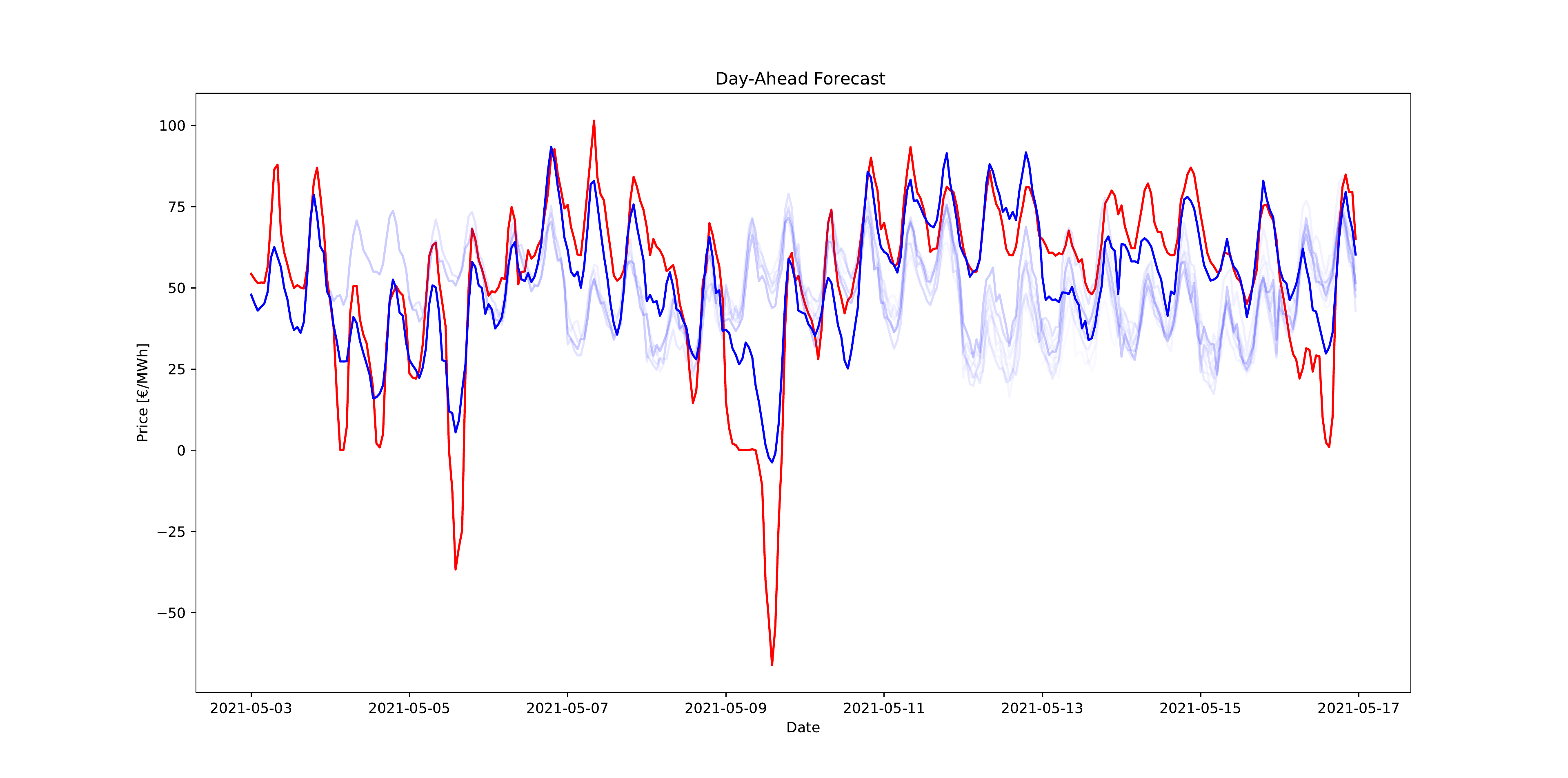}
    \caption{Realized day-ahead prices (red) and daily day-ahead forecasts (blue) for a two week time period at the start of May 2021.}
    \label{fig: day-ahead example}
\end{figure}

For the price forecasts, we use the forecasting model from Wagner et al. (cf. \cite{ramentol2020short}). The model is based on a simple dense neural network (DNN) with an embedding layer to process calendar data. The underlying data is the date and time of the price to be predicted and the corresponding forecast of wind and solar input. The DNN is trained daily on the five years preceding the two weeks we want to forecast, with the above data as input and the day-ahead prices as output. We obtain the data from the ENTSO-E transparency platform (cf. \cite{entsoe}). As the data for wind and solar are only available for the next day, but we are forecasting for two weeks, we make two separate forecasts. First, we forecast the next 24 hours using the ENTSO-E renewables data. Then we forecast the next 13 days without this data using a DNN model trained only on calendar information. This results in a day-ahead price forecast for the next two weeks. For an example of the forecasted prices see Figure~\ref{fig: day-ahead example}. Note that every forecast gets more opaque over time. So the solid blue line represents the first 24 hours of each forecast.

\section{The cross-market optimization problem}\label{sec: cm problem}
In this section we give a detailed description of our cross-market optimization problem. Note that we are considering a rolling horizon approach as described in Section~\ref{subsec: rolling horizon}. This requires us to solve an optimization problem every day. We now describe this optimization problem for the two-market setting, where we aim for a combined bidding strategy on the German balancing markets and the German day-ahead market. The solution obtained is not only a distribution of flexibility between these markets, but also a precise bidding strategy for each of the markets. 

The variables are introduced in the following. A complete formulation of the optimization problem can be found in Appendix~\ref{app: problem}. We begin by introducing some general notation. We denote the number of hours considered for cross-market trading by $T$ and the number of hours considered for the total look-ahead period by~$T^*$. This is done to allow a longer view of the day-ahead market and a penalty term, which we will introduce later. Notation regarding products sets was already introduced in Section \ref{sec: markets}.

Furthermore, the expected profit on each market is denoted by $f^\RL$ (or $f^\RA$, $f^\DA$). We introduce decision variables~$m^\text{i}$ for $i\in\{\DA, \RM\}$, which denotes the flexibility allocated to the respective market. For better readibility, we introduce variables~$m^\text{i}$ for $i\in\{\RL, \RA\}$ as well and set
\begin{align}
    m^{\RM} = m^{\RL} = m^{\RA}.
\end{align}
Furthermore, we add a \emph{penalty term} regarding battery usage, which we denote by $f^\TC$. 

\paragraph{The objective function}
The objective function of our cross-market model is given as the sum of the expected profits in the single markets minus the penalty term, i.e.
\begin{align}\label{eq: cm_opt}
    \max \quad f^\RL + f^\RA + f^\DA - f^\TC. 
\end{align}

In the following, we describe the expected profit of each of the markets in more detail and formulate corresponding market constraints, as well as linking constraints among the decision variables of all markets.  

\paragraph{Reserve capacity market constraints}
As described in Section~\ref{subsubsec: reserve capacity}, for a given day we compute $N^{\RL}$~bidding scenarios that are given by price levels~$p_i^\RL$ with $p_i^\RL \leq p_{i+1}^\RL$ for all $i=1,\dots, N^{\RL}-1$ and the corresponding acceptance probabilities~$q_i^\RL$. Given these, the \emph{expected profit} of product~$k$ is defined as
\begin{align}
    f_k^\RL \coloneqq \sum_{i=1}^{N^{\RL}} q_{k, i}^{\RL} \left(\sum_{j=1}^i m_{k, j}^{\RL}\cdot p_j^\RL\right),
\end{align}
where $m_{k, i}^{\RL}$ denotes the offered volume at price level~$p_i^\RL$ and for product~$k$. The total available volume for product~$k$ is bounded by the distributed flexibility~$m_k^\RL$ which again is bounded by the allocated flexibility~$m^\RL$. Thus, the constraints
\begin{align}
    \sum_{i=1}^{N^{\RL}} m_{k, i}^{\RL} &\leq m^\RL  & \forall k\in \mathcal{K} \\
    m_k^\RL &\leq m^\RL & \forall k\in \mathcal{K}
\end{align}
are added. 
Recall, the \emph{total expected profit} of the \RL~market is denoted by the variable~$f^\RL$ and is defined as the sum of expected profits over all products in $\mathcal{K}$, i.e. 
\begin{align}
    f^\RL = \sum_{k \in \mathcal{K}} f_k^\RL.
\end{align}

\paragraph{Reserve energy market constraints}
In Section~\ref{subsubsec: reserve energy} we described how to compute the \emph{expected activation duration}~$\alpha_{k, i}$ for a given price level~$p_i$ and a product~$k$. We consider $N^{\RA}$ different price levels. As said before, the considered price levels of the \RA~market do not need to be the same as the considered price levels on the \RL~market. Now, given these values, the \emph{expected profit on the \RL~market} of product~$k$ is defined as 
\begin{align}
    f^\RA_k \coloneqq \sum_{i=1}^{N^{\RA}} p_i^\RA\cdot m_{k, i}^{\RA}\cdot \alpha_{k, i},
\end{align}
and the \emph{total expected profit of the \RL~market} is the sum of expected profits over all products, i.e. we add the constraint
\begin{align}
f^\RA = \sum_{k\in \mathcal{K}} f^\RA_k. 
\end{align}

\paragraph{Day-ahead market constraints}
In order to compute the expected day-ahead profit, recall that the total flexibility distributed to the day-ahead market is denoted by $m^\DA$. Since we consider only products that are traded at the auction for hourly products, we introduce a decision variable~$m_t^\DA$ for every hour~$t=0,\dots,T^*$. Furthermore, we introduce a decision variable~$m_k^\DA$ for each 4-hour block~$I_k\in \mathcal{I}$. We link them by the constraints
\begin{align}
    -m_k^\DA &\leq m_t^\DA \leq m_k^\DA &\forall I_k\in \mathcal{I}, t\in I_k, \\
    -m^\DA &\leq m_k^\DA \leq m^\DA  &\forall k\in \mathcal{K}.
\end{align}

Consequently, $m_t^\DA$ denotes the amount of MW that is offered on the day-ahead market during hour $t$ (where a negative amount translates to a buying bid). The constraints ensure that this amount cannot be higher than the maximally allocated flexibility for this market.

Given hourly price forecasts~$[p_t]$ for $t = 0,\ldots,T^*$, we define the expected profit on the day-ahead market as
\begin{align}
    f^\DA & \coloneqq \frac{T}{T^*}\sum_{t=0}^{T^*} p_t \cdot m_t^\DA, 
\end{align}
where $p_t$ denotes the price for the hourly product at hour $t$. Note that we downscale the profit made on the day-ahead market to the considered time period for the cross-market trading in order to obtain profits that are comparable to the balancing market profits. Recall Section~\ref{subsec: da market} on how the forecasts are computed. 

\paragraph{Battery capacity constraints}
We denote the minimum (respectively, maximum) battery capacity by $\tau_{\min}$ (respectively, $\tau_{\max}$). At each point in time~$t$, the battery level has to lie in the interval~$[\tau_{\min}, \tau_{\max}]$. The initial battery level is denoted by $\tau_0$. In total, we obtain the following battery capacity constraints. 
\begin{align}
    \tau_{\min} \leq \tau_{0}  + \sum_{t'=0}^{\min\{t,T\}} m_{t'}^\RL + \sum_{t'=0}^{t} m_{t'}^\DA &\leq \tau_{\max} \qquad \forall t = 0,\dots,T^*, \label{cons: bat cap 2}
\end{align}
where $T$ denotes the hours of the first day whereas $T^*$ denotes the number of hours of the whole considered period. Note that these constraints ensure that the battery level remains in the total capacity interval at every hour and link the flexibility given to each market in one constraint. Furthermore, it should be noted that the distributed flexibility to the balancing markets is considered in this constraint by~$m_t^\RL$. Again, the usage of $m_t^\RL$ represents a conservative estimate, as this actually equals the maximum amount of energy that could possibly be discharged from the battery through the \RA~market.

For a short-term consideration, we could leave the objective function as $\max f^\RL + f^\RA + f^\DA$. However, in a medium or long term setting, this will often result in a short-sighted exploitation of the battery capacity. This is especially undesirable if the production efficiency (as in our application use case) suffers significantly from a deviation of the battery from its central position.

In order to prevent discharging strategies, the next two paragraphs introduce constraints that counteract this, namely a penalty for high absolute battery deviations and the constraint to obtain a balanced battery level at the end of the look-ahead period. Note that we cannot control future battery levels precisely, as the future state depends on the realisations of bids on the balancing markets and is also being overturned by interim scheduling updates. Hence we can only make assumptions about the expected battery states.

\paragraph{Penalty term constraints}
A high efficiency cost parameter will prevent trading strategies that results in very high or very low battery levels since any deviation from the nominal baseline is penalized. For this, we denote the (expected) battery level at time point~$t$ by $\tau_t$ and define it as 

\begin{align}
    \tau_t \coloneqq \sum_{I_k\in \mathcal{I} : a_k\leq t} \left( \sum_{i=1}^{N^{\RA}} m_{k, i}^{\RA}\cdot \alpha_{k, i} \right) + \sum_{t'=0}^t m_{t'}^{\DA} \label{cons:penaltyterm1}
\end{align}

Due to constraints~\eqref{cons: bat cap 2} on the battery capacity, the expected battery level always lies in the interval~$[\tau_{\min}, \tau_{\max}]$. We now define the technical cost penalty term as
\begin{align}
    f^\TC = \frac{T}{T^*}\sum_{t=0}^{T^*} c_p |\tau_t|, \label{cons:penaltyterm2}
\end{align}
where $c_p$ denotes the \emph{penalty coefficient}. We will later evaluate strategies with respect to different penalty coefficients. Note that a higher penalty coefficient yields strategies that balance the battery in shorter periods. Here, we considered a linear formulation in order to obtain a linear program. However, one might consider a quadratic penalty term in order to penalize higher (or lower) battery levels more. Furthermore, note that we again use a scaled approach similar to the day-ahead profit function. 

\paragraph{Battery balance constraints}
As pointed out before, we compute day-ahead strategies over the longer period~$T^*$. In order to prevent selling all remaining energy in the battery at the end of $T^*$ (as this would maximize the profit), we require that the expected battery level at the last time point is zero, i.e.
\begin{align}
    \tau_{T^*} = 0.
\end{align}

\paragraph{Flexibility constraints}
It is necessary to globally restrict the flexibility distributed to the different markets such that it cannot exceed the total flexibility available. This is ensured by the constraint
\begin{align}
    m^\DA + m^\RM \leq m,
\end{align}
where $m$ denotes the total available flexibility that represents the maximum charging or discharging amount of electricity.

This concludes our discussion of the constraints. A full formulation of the problem can be found in Appendix~\ref{app: problem}. 

\section{Numerical case study}\label{sec: numerical case study}
In this section, we give a brief description of our case study setting. At TRIMET Aluminium SE, the aluminium production process was flexibilized such that the process now functions as a \emph{virtual} battery (cf. \cite{dussel2019transformation, Depree2016}). 

The production of aluminium typically involves electrolytic reduction, which is a very energy intensive process. In this process, aluminium is extracted from alumina by passing an electric current through a mixture of alumina and cryolite held in reduction cells. The process has traditionally been run at a constant high current, as any change in current can potentially lead to cell damage. More details about the aluminium production process can be found in~\cite{Depree2016, dussel2019transformation}. 

TRIMET Aluminium SE, based in Essen, Germany, has developed a pilot set of reduction cells equipped with heat exchangers. These devices enable a more flexible energy schedule, as they allow the energy input to be adjusted within certain time and energy constraints while keeping the risk of cell damage reasonably low. This set-up can act as a \emph{virtual battery}, as a surplus or shortage of energy can be stored in the cells and balanced within a given time (cf.~\cite{Depree2016, schafer2019model}). This allows the company to take advantage of price spreads on different energy markets. A more detailed explanation of the \emph{virtual battery} can be found in \cite{Depree2016, dussel2019transformation}. For the proposed model, the virtual battery has a capacity of 1000~MWh and a total flexibility of 10~MW per hour, i.e. we could sell or buy energy for 100 consecutive hours. Note that the parameter choices of the virtual battery are inspired by the setting at TRIMET Aluminium SE, but do not reflect their real parameter values.

As described in constraints~\eqref{cons:penaltyterm1} and \eqref{cons:penaltyterm2}, we define different bidding strategies by two parameters, the flexibility distribution~$(m^{\DA}, m^{\RM})$ between day-ahead market and balancing market and the penalty coefficient~$c_p$. If these two parameters can be adjusted between different runs during the rolling horizon computation, we simply denote this with \texttt{free}. For example, the strategy that has a fixed distribution of 4~MW to the day-ahead market and 6~MW to the balancing markets with a fixed penalty coefficient of 0 is denoted ((4,6) $|$ 0), whereas the strategy that can adjust both the flexibility distribution and the penalty coefficient over time is denoted (\texttt{free} $|$ \texttt{free}). 

\begin{figure}[ht]
    \centering
    \includegraphics[width=12cm]{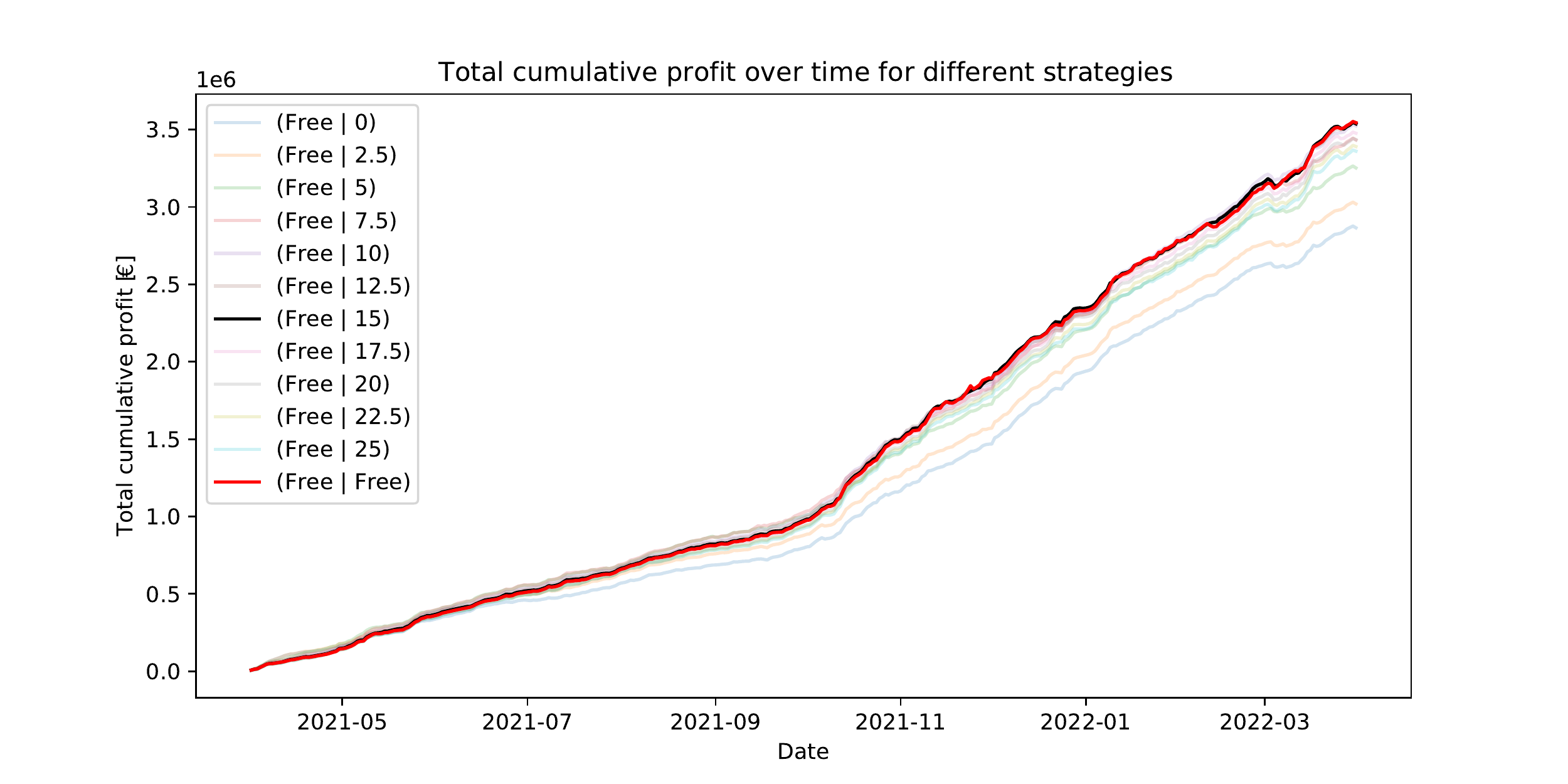}
    \caption{Cumulative realized profits for different (\texttt{free} $|$ $\ast$)--strategies. The (\texttt{free} $|$ \texttt{free})--strategy is shown in red and obtains a profit of approximately 3.551.222~\euro. The (\texttt{free} $|$ 15)--strategy, which obtains the highest profit of all free strategies with a fixed penalty coefficient, is shown in black. All other strategies are colored with a lower opacity.}
    \label{fig: tm total profits}
\end{figure}

A high penalty coefficient keeps the battery utilization within a narrower range, i.e. by penalizing very large deviations from a balanced battery level. In the following, we will not only examine strategies with a fixed coefficient, but also a strategy that adjusts this coefficient on a daily basis according to the current market situation. To do this, we clustered various price and volatility metrics on the day-ahead and balancing markets and now choose the coefficient according to the cluster it belongs to. Our strategy is to choose a low coefficient when volatility tends to be low, in order to have the possibility of realising profits over a longer period of time, corresponding to a higher utilization of the battery. 

\begin{figure}[ht]
    \centering
    \includegraphics[width=12cm]{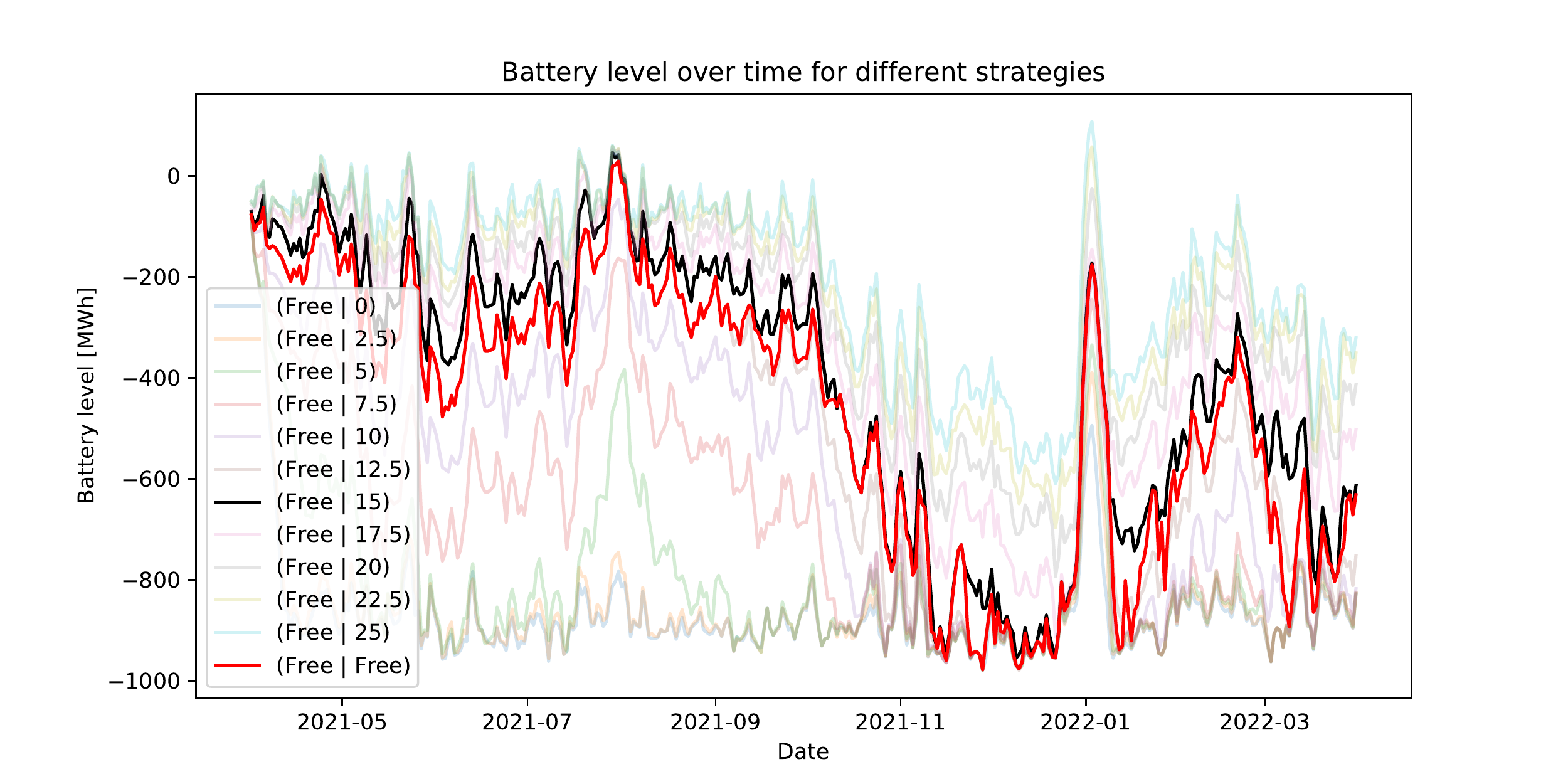}
    \caption{The battery level over time for different (\texttt{free} $|$ $\ast$)--strategies. The best (in terms of total cumulative profit) strategy with fixed penalty coefficient is shown in black, whereas the (\texttt{free} $|$ \texttt{free})--strategy is shown in red. Overall, we observe a trend to discharge the battery towards the end of 2021. The high spike around new years eve is due to an unusual high demand of negative \RA~energy.}
    \label{fig:batterylevel}
\end{figure}

\begin{figure}[ht]
    \centering
    \includegraphics[width=10cm]{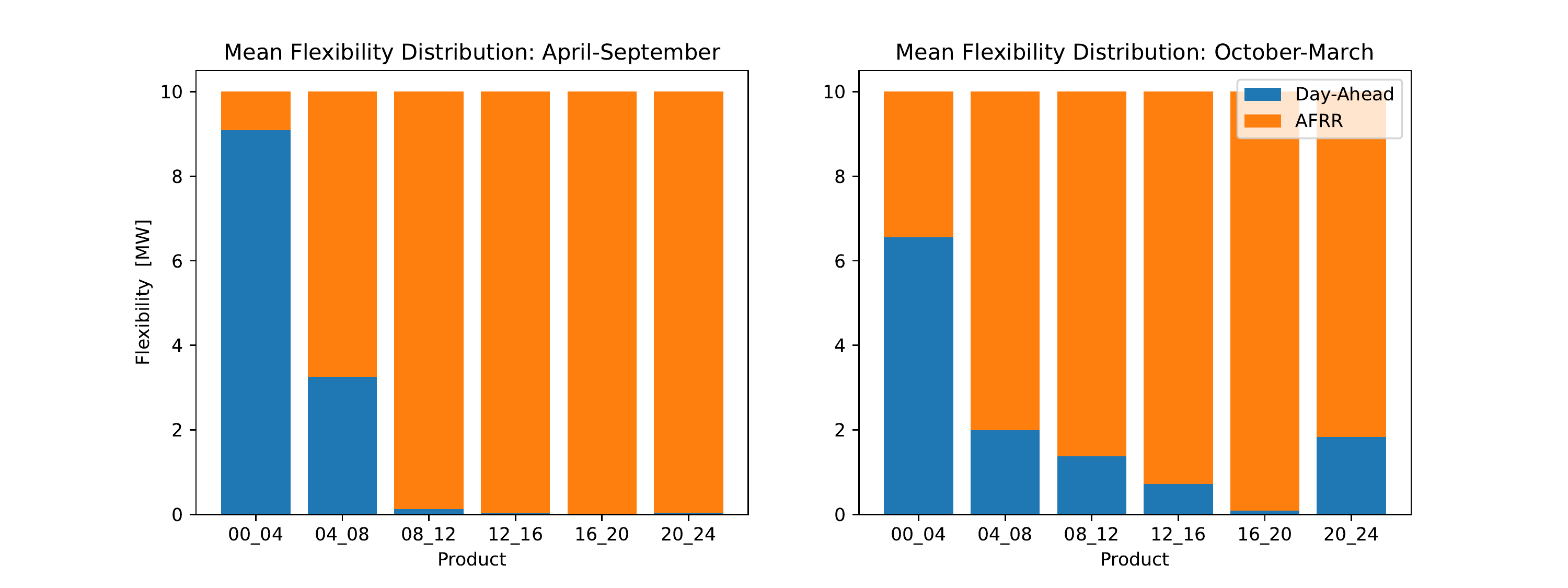}
    \caption{Flexibility distributed to the reserve markets and day-ahead market for the (\texttt{free} $|$ \texttt{free})--strategy split between the first, i.e. from 04/2021 to 11/2021, and second, i.e. from 12/2021 to 03/2021, of the evaluation period. We observed a steeper cumulative profit curve during the second period which indicates a structural change of market prices. The figure shows that this also changes the bidding strategy.}
    \label{fig: mean split distribution periods}
\end{figure}

We observe very similar cumulative profits for all considered strategies as shown in Figure~\ref{fig: tm total profits}. What is striking is that, depending on the obtained profit, the period can be divided into two sub-periods. For the whole period, the (\texttt{free} $|$ \texttt{free})--strategy (colored red) achieves the highest total profit of 3.551.222~\euro, although it is only about 1\% higher than the (\texttt{free} $|$ 15)--strategy (colored black). An overview over all profits obtained by the different bidding strategies can be found in Table~\ref{table:overview profits} in the Appendix. The (\texttt{free} $|$ \texttt{free})--strategy shows minor differences from the (\texttt{free} $|$ 15)--strategy, tending to keep the battery at a slightly lower level. Around the turn of the year, we see a small drop in profit, indicating that we are buying energy over several consecutive days. We can see this more clearly in Figure~\ref{fig:batterylevel}, which shows battery levels over time. This effect is more obvious at lower penalty coefficients. It is interesting to note the steeper curve starting in October 2021, indicating increased profit opportunities. In Figure~\ref{fig: mean split distribution periods}, we see that the average allocation of flexibility in the first period is similar, but not identical to that in the second period. The higher profits in the second period go hand in hand with higher price-bids in the positive SRE market. While the average bid in the first period is 130~\euro/MWh, it then more than doubles to 287~\euro/MWh in the second period.  In comparison to fixed strategies, i.e. strategies that market the same flexibility distribution every day, we see that all strategies that allow a daily adjusted split outperforms these fixed strategies, as shown in Figure~\ref{fig: fixed strategies} and Table~\ref{table:overview profits}.

Figure~\ref{fig:batterylevel} depicts the effects of the different penalty terms on the trading decisions. Starting with an battery level of 0~MWh, the strategies with a low penalty term instantly start to heavily discharge the battery. Although that yields a high profit at first, being close to the minimum limits the possibilities of trading. With an increasing penalty term, the discharging effect reduces. For the highest penalties, the battery level stays close to zero - at least for the first half. Thereafter, the electricity prices start to increase strongly and thus also profit margins. None of the penalties suffices to prevent the battery from discharging in all the strategies. At the end of 2021, even the (\texttt{free} $|$ 25)--strategy reaches a minimum of around -600~MWh. The rest of the figure depicts the more volatile prices in the beginning of 2022. The prices in the first two months of 2022 start to relax, accompanied with an overall increase of battery levels for the strategies. At the end of February, the reverse is true and the battery levels start to fall again. The jump around new years is due to an unusual high demand on the negative \RA~market, see e.g. Figure~\ref{fig: tm retrieval times} in the appendix. The best performing strategies -- (\texttt{free} $|$ \texttt{free}) and  (\texttt{free} $|$ 15) -- are located neither at the top nor at the bottom of the range of realized battery levels. This hints at the conclusion that a well calibrated penalty term should be high enough so that it does not allow for complete discharge, while being low enough such that trading is not restricted too much.

\begin{figure}[ht]
    \centering
    \includegraphics[width=12cm]{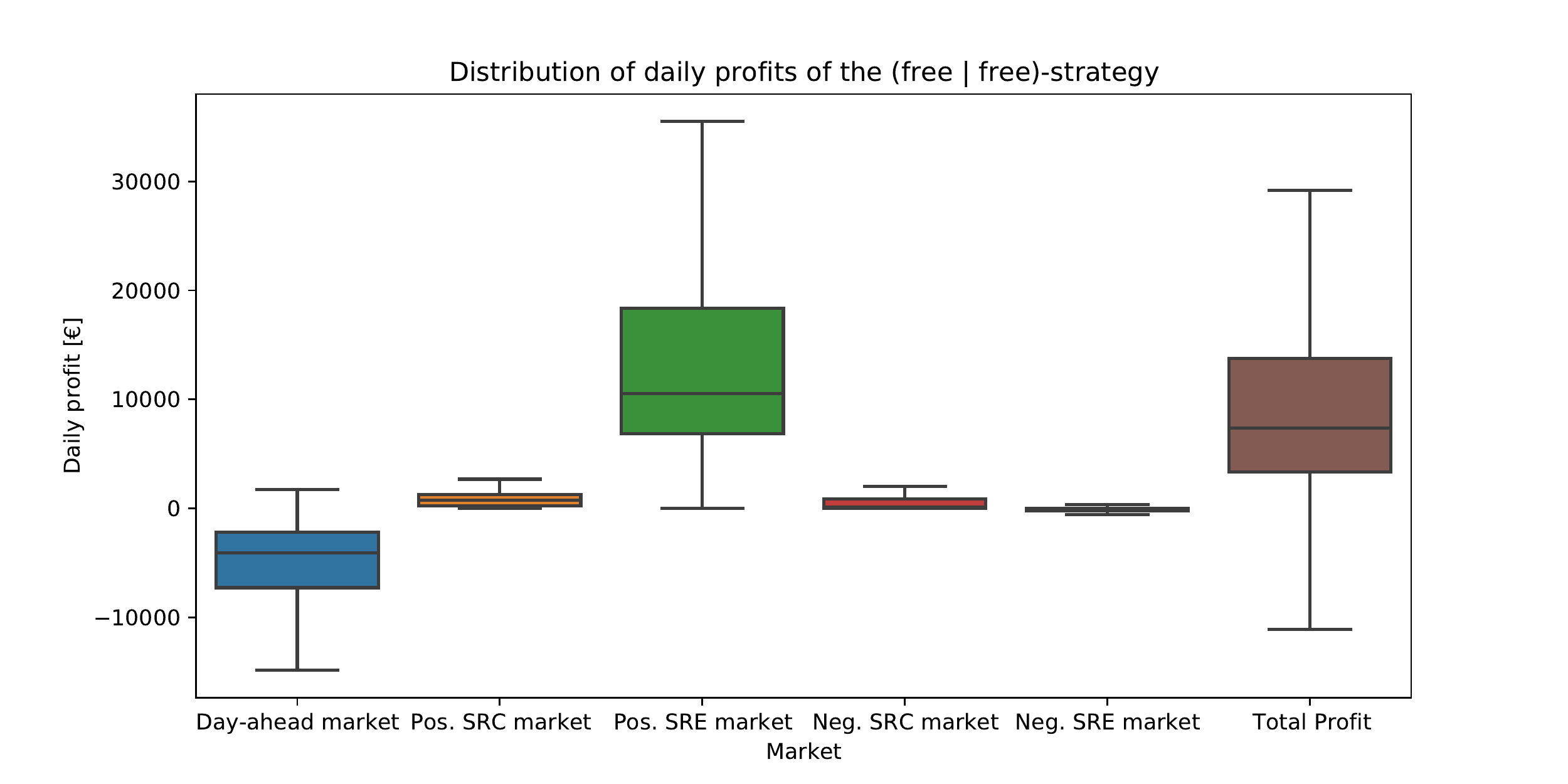}
    \caption{Distribution of the realized daily profits of the (\texttt{free} $|$ \texttt{free})-strategy. Note that while the day-ahead market profits are mainly negative, i.e. energy is bought, the main profit contribution comes from the positive \RA~market. Since we allow both positive and negative price bids on the negative \RA~market, the daily profit can be both negative and positive, unlike the other balancing markets where the profits are non-negative.}
    \label{fig: tm daily total profits}
\end{figure}

In order to gain an intuition about realized profits of the best-performing (\texttt{free} $|$ \texttt{free})--strategy, its daily profits are displayed in a box-plot in Figure~\ref{fig: tm daily total profits}. Several key values are presented there: For one, each box represents the daily profits that were in between the 25\% quantile and the 75\% quantile for the corresponding market. The black line in the middle of the blue box indicates the 50\% quantile, which represents the value that lies in the middle of all observed and ordered daily profits. Consequently, each box contains 50\% (which equals 75\% - 25\%) of all observed profits over the considered time horizon for the respective market. The vertical lines above and below each box mark the distance of the 25\% quantile (respectively the 75\% quantile) to 1.5 times the distance between the 25\% quantile and the 75\% quantile. The so-called `whiskers' at their end are then either determined by this length or by the furthest data point available in the corresponding direction. In Figure~\ref{fig: tm daily total profits}, it becomes clear that the day-ahead market is mainly used to buy electricity, because most of its profits are negative. In comparison to that, the positive \RA~market is the market yielding the highest profits and is per definition used to sell electricity, while the negative \RL~market is used to refill the battery. Both \RL~markets, as well as the negative \RA~market, do not yield profits that are of much significance, nor do they show much variability in the attained values. Consequently, the total profit's volatility is mainly influenced by day-ahead and positive \RA~market volatility. Furthermore, as the quantiles of the total daily profit are well above zero, it must hold true that the overall usage of positive \RA~is generally high, as it is the only market contributing strongly to positive daily profits.

\begin{figure}[ht]
    \centering
    \includegraphics[width=12cm]{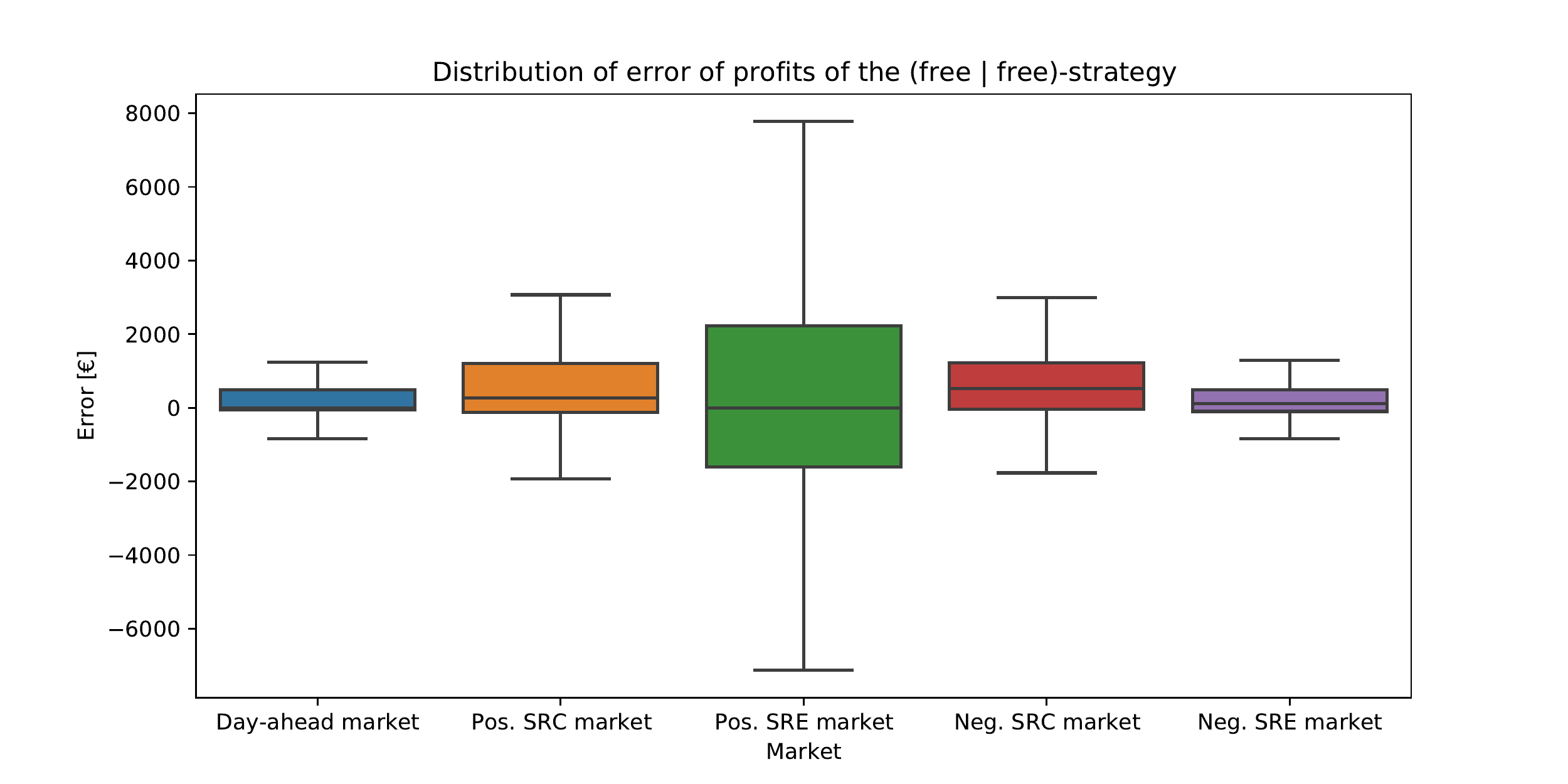}
    \caption{Difference of expected daily profit minus realized daily profit for different markets. We see that the error is highest for the positive \RA~market, which is also the market with the highest share of the total profit.}
    \label{fig: error daily profits}
\end{figure}

When speaking of daily realized profits, one also has to consider the expected daily profits as shown in Figure~\ref{fig: error daily profits}. Albeit being the market with the highest daily profits, the positive \RA~market is also the market with the highest difference between expected and realized profit. 
We see that the median difference is almost zero, which means the forecast for the \RA~market is roughly equally often too big as too small. 
Since the obtained profit of the other balancing markets is significantly lower, the difference does not have too much impact. We see that the day-ahead forecasts perform really well with a small shift towards too optimistic expected profits. Overall, our trading strategy could be improved by using better forecasts, in particular for the \RA~markets.

\begin{figure}[ht]
    \centering
    \includegraphics[width=12cm]{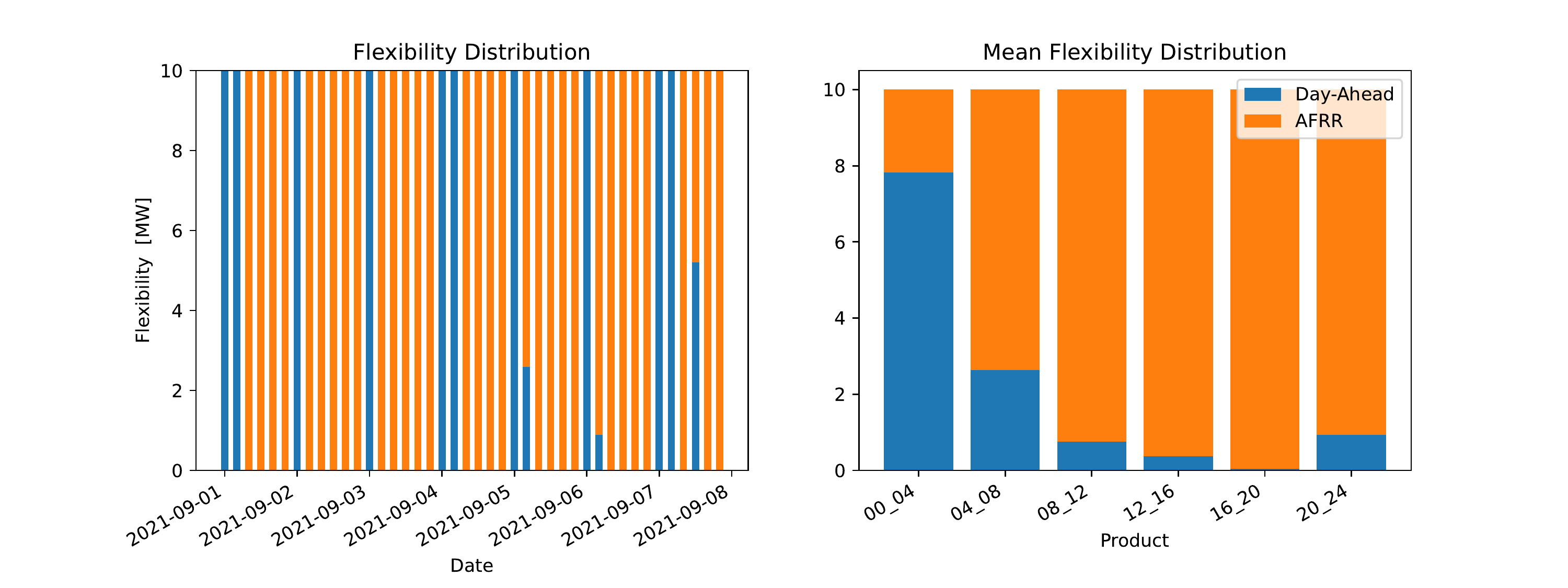}
    \caption{Flexibility distributed to the \RM~and day-ahead markets for the  (\texttt{free} $|$ \texttt{free})--strategy. The left panel shows the first seven days of 09/2021, the right panel shows the mean distributed flexibility per product over the complete evaluation period.}
    \label{fig: tm flexibility distribution}
\end{figure}

We show in Figure~\ref{fig: tm flexibility distribution} the daily split of flexibility of the (\texttt{free} $|$ \texttt{free})--strategy. In the left panel the distribution per 4h-block is shown for the first seven days of September 2021. Note that the split is not always trivial, i.e. the strategy distributes flexibility to both markets in the same 4h-block. This is due to the fact that the strategy does not solely allocate the flexibility to the market with the highest expected profit, but also takes the penalty term into account. The right panel of Figure~\ref{fig: tm flexibility distribution} shows the mean distribution of flexibility over the whole period considered. It is striking, that the balancing markets are strongly preferred over the day-ahead market. However, in the first four hours of the day, the strategy allocates the available flexibility mostly to the day-ahead market. From 4~a.m. to 8~a.m. a quarter of the flexibility is still offered on the day-ahead market. This fits with the result that the day-ahead market is mainly used for buying, as these hours usually exhibit the lowest prices. For all other 4h-blocks, the balancing markets are strongly favored.

Overall, we see that all strategies that adjust the flexibility distributed to the different markets outperform all strategies with a fixed flexibility split. By also adjusting the penalty coefficient on a daily basis, the profit can be further increased. However, our adjustment provides only a slight increase, but the potential of a daily adjusted penalty coefficient becomes obvious. 

\section{Conclusion}\label{sec: outlook}
In this article, we propose a solution to the optimization problem of distributing flexibility between the day-ahead market and the \RM~market in the setting of a virtual battery. We consider the case where a decision on the flexibility distribution has to be made before the first market closes in order to guarantee real-life applicability of our solution. In this setting, we only solve one optimization problem only per day which takes both markets into account, in contrast to other approaches that solve the problem in a multistage setting. The resulting flexibility distribution cannot be changed afterwards. We present different trading strategies for these markets that are reached by optimizing over various flexibility distributions between the markets plus a penalty coefficient.

A numerical study over the period from April 2021 to March 2022 shows that a daily adjusted split of the available flexibility yields significantly higher profits than a strategy that tries to market the same split of flexibility every day on the day-ahead and the secondary balancing markets. Furthermore, we introduce a penalty coefficient that determines how much deviations from the base battery line are penalized. 
Our proposed (\texttt{free} $|$ \texttt{free})--strategy, which can adjust both the flexibility and the penalty coefficient on a daily basis, outperforms all other strategies considered in terms of total profits. By adjusting the penalty coefficient between daily optimization runs, we obtain a desired property of the resulting plan, i.e. maximum profit. While we have presented a penalty approach that depends on the volatility measures of the considered markets, one can extend this by considering, for example, the best choices of the coefficient over the last few days, the current battery level, and other external factors such as solar and wind forecasts. This is part of future research. Furthermore, finding an appropriate adjustment of the penalty coefficient is a crucial task that also allows further improvement.

Another market that offers flexibility is the continuous intraday market. Therefore, a further possible extension of the model is the inclusion of this market. The intraday market differs from the others in that buyers and sellers of electricity meet directly. This means that they place bids which are collected in an order book and when two orders for a product match, i.e. the selling price is lower than or equal to the buying price, a trade is executed. Our proposed method of incorporating the intraday market is through an algorithm that tracks the order book and trades arbitrage between different products. The profit generated by this algorithm is then used as a time series to forecast profits for cross-market optimization. Difficulties arise in calibrating the algorithm to produce results that reflect possible outcomes in a real-world scenario. In addition, the intraday market is much more volatile for participants than the other markets, making quality forecasting more difficult.

\section*{Acknowledgement}
This work has partially been supported by the German Federal Ministry for Economics and Climate Action in grant 01186724/1 (FlexEuro: Wirtschaftliche Optimierung flexibler stromintensiver Industrieprozesse). 

\bibliographystyle{plain}
\bibliography{ref}

\appendix

\section{Two market problem formulation}\label{app: problem}
Here we give a full formulation of the two market optimization problem. 
{\small 
\begin{align*}\label{eq: cm_opt appendix}
    \max \quad &  f^\RL + f^\RA + f^\DA - f^\TC \\
    \text{s.t.} \quad & f^\RL = \sum_{k \in \mathcal{K}} \sum_{i=1}^N q_{k, i}^{\RL} \left(\sum_{j=1}^i m_{k, j}^{\RL}\cdot p_j^\RL\right) \\
    & f^\RA = \sum_{k\in K} \sum_{i=1}^{N^{\RA}} p_i^\RA\cdot m_{k, i}^{\RA}\cdot \alpha_{k, i} \\
    & f^\DA = \frac{T}{T^*}\sum_{t=0}^{T^*} p_t \cdot m_t^\DA \\
    & f^\TC = \frac{T}{T^*}\sum_{t=0}^{T^*} c_p |\tau_t| \\
    & \tau_{\min} \leq \tau_{0} + \sum_{t'=0}^{\min\{t,T\}} m_{t'}^\RL + \sum_{t'=0}^{t} m_{t'}^\DA \leq \tau_{\max} & \forall t' = 0,\dots,T^* \\
    & \tau_t = \sum_{I_k\in \mathcal{I} : a_k\leq t} \left( \sum_{i=1}^{N^{\RA}} m_{k, i}^{\RA}\cdot \alpha_{k, i} \right) + \sum_{t'=0}^t m_{t'}^{\DA} & \forall t=0,\dots, T^*\\
    & \tau_{T^*} = 0 \\
    & m^{\RM} = m^{\RL} = m^{\RA} \\ 
    & m^\DA + m^\RM \leq m \\
    & -m_k^\DA \leq m_t^\DA \leq m_k^\DA &\forall k\in \mathcal{K}, t\in I_k \\
    & -m^\DA \leq m_k^\DA \leq m^\DA  &\forall k\in \mathcal{K} \\
    & \sum_{i=1}^{N^{\RL}} m_i^{\RL, k} \leq m^\RL  & \forall k\in \mathcal{K} \\
    & m_k^\RL \leq m^\RL & \forall k\in \mathcal{K}
\end{align*}
}

\section{Retrieval times on the \RA~markets}
\begin{figure}[H]
    \begin{minipage}{.5\textwidth}
        \centering
        \includegraphics[width=7cm]{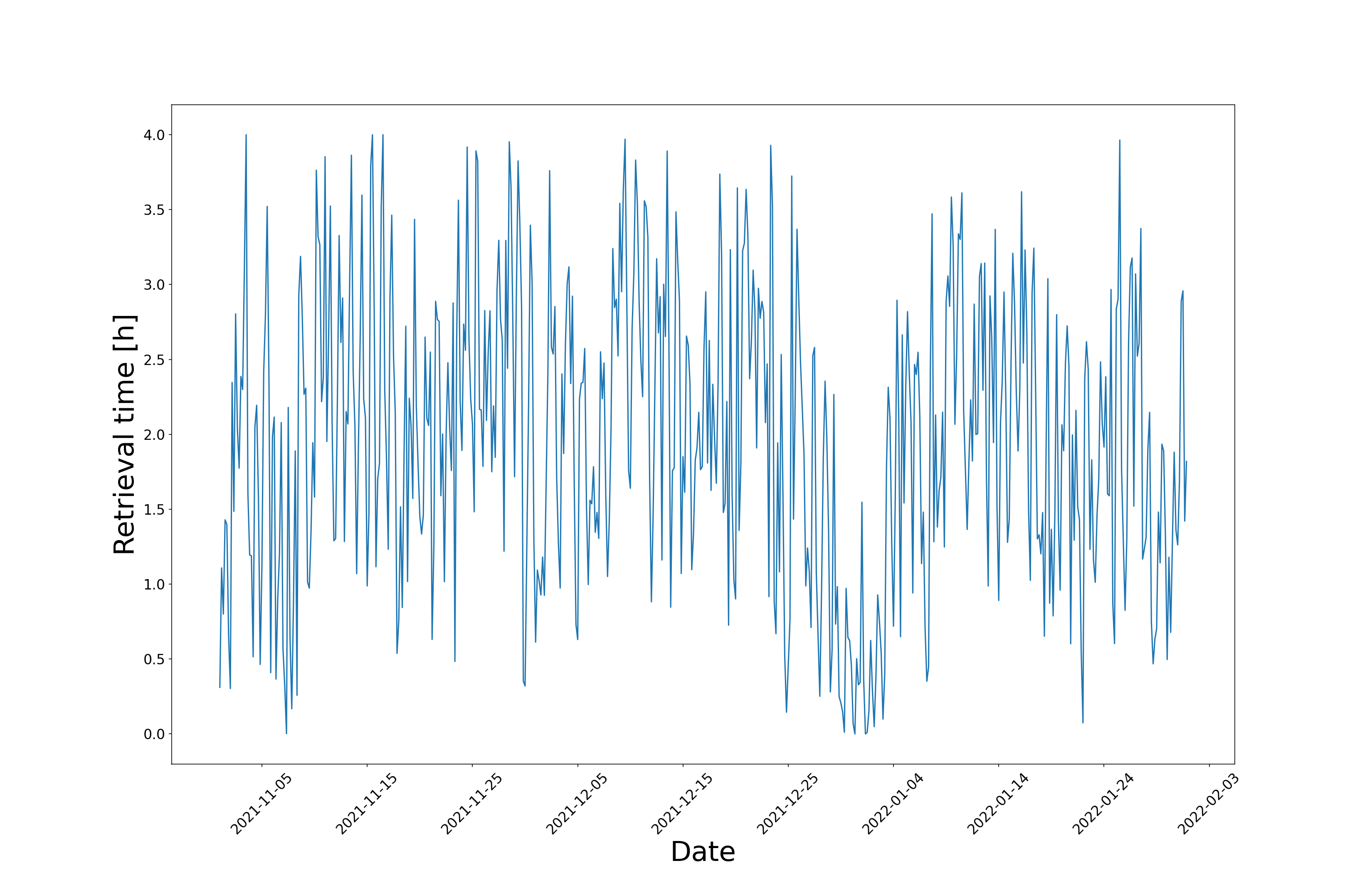}
    \end{minipage}%
    \begin{minipage}{.5\textwidth}
        \centering
        \includegraphics[width=7cm]{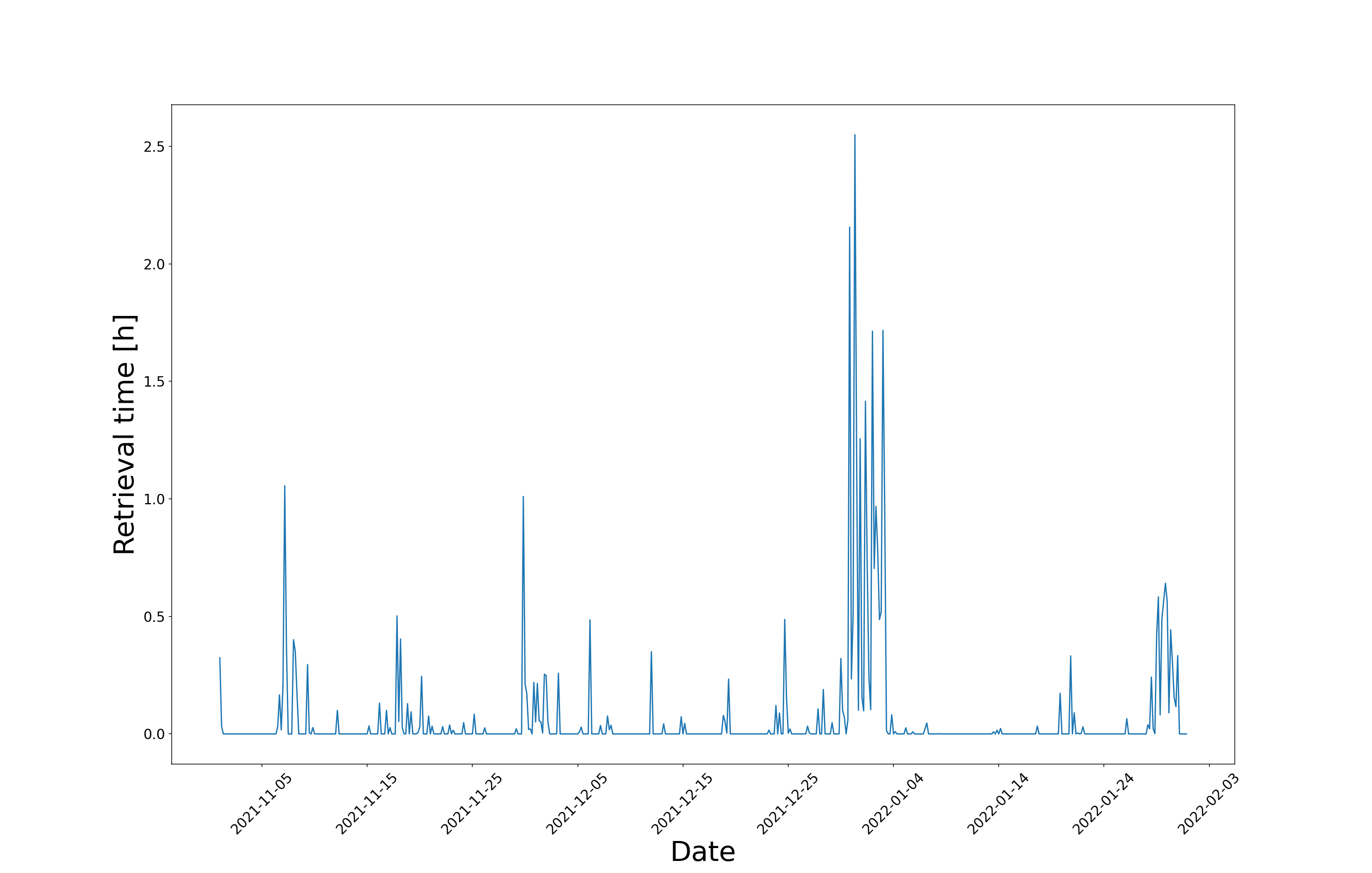}
        
    \end{minipage}%
    \caption{Retrieval times of the positive reserve energy market for price level of 95~\euro/MWh. We observe an oscillating pattern of retrieval times in the time period from Nov 2021 - Jan 2022, where a significantly low demand during Christmas/New Year season in 2021 becomes apparent. Here the x-axis is based on a 4h-scaling as the product structure of the market in that time period would suggest.} \label{fig: tm retrieval times}
\end{figure}
\section{Total profits of different bidding strategies}
\begin{table}[H]
\centering
{\small 
\begin{tabular}{cccl}
\hline\\
Flexibility DA [MW]& Flexibility aFFR [MW]& Penalty coefficient &      Profit [\euro] \\
\hline\\
          free &             free &                free &  3.551.222 \\
          free &             free &                  15.0 &  3.537.889 \\
          free &             free &                12.5 &  3.536.887 \\
          free &             free &                  10.0 &  3.523.781\\
          free &             free &                17.5 &  3.476.709\\
          free &             free &                  20.0 &  3.435.494 \\
          free &             free &                 \phantom{0}7.5 &  3.433.082 \\
          free &             free &                22.5 &  3.390.784\\
          free &             free &                  25.0 &  3.360.793 \\
          free &             free &                   \phantom{0}5.0 &  3.252.534\\
          free &             free &                 \phantom{0}2.5 &  3.020.029 \\
          free &             free &                   \phantom{0}0.0 &  2.865.539 \\
           4.0 &              6.0 &                   - &  2.500.469\\
           5.0 &              5.0 &                   - &  2.452.831\\
           3.0 &              7.0 &                   - &  2.439.196 \\
           6.0 &              4.0 &                   - &  2.348.071 \\
           7.0 &              3.0 &                   - &  2.244.507 \\
           8.0 &              2.0 &                   - &  2.104.400\\
           2.0 &              8.0 &                   - &  2.010.676 \\
           9.0 &              1.0 &                   - &  1.947.897 \\
          10.0 &              0.0 &                   - &  1.775.500 \\
           1.0 &              9.0 &                   - &  1.129.863 \\
           0.0 &             10.0 &                   - &   \phantom{0.}115.723\\
\hline\\
\end{tabular}
}

\caption{Overview over total obtained profits at the end of the considered period for the different bidding strategies, ordered by the profit in ascending order.}
\label{table:overview profits}
\end{table}

\section{Evaluation of (fixed $|$ 0)-strategies for different flexibility splits}

\begin{figure}[H]
    \centering
    \includegraphics[width=10cm]{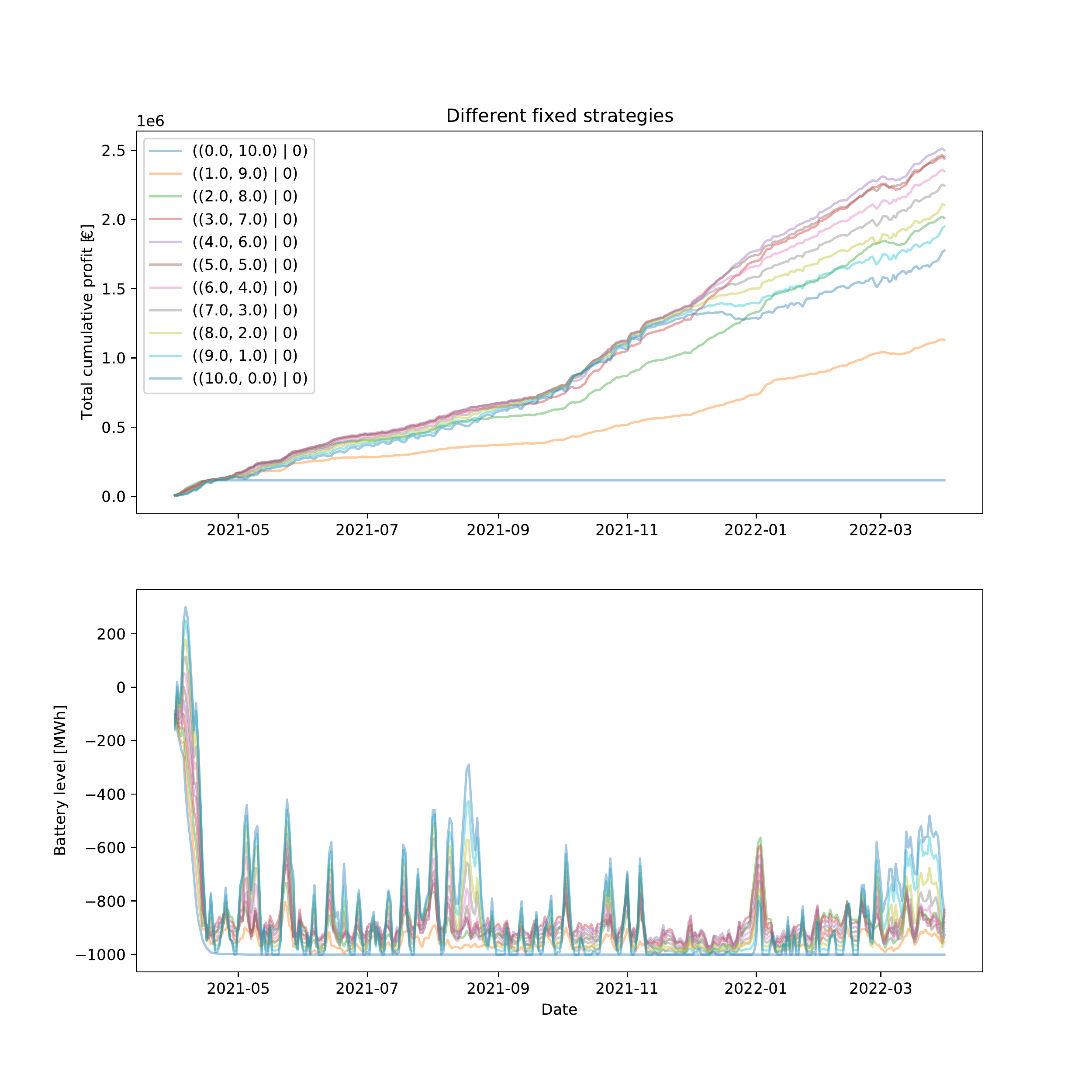}
    \caption{Evaluation for different (fixed $|$ 0)-strategies. In terms of total cumulative profit, the fixed split of 4~MW on the day-ahead market and 6~MW on the balancing markets outperforms all other strategies. Note that strategies with a high distributed amount of flexibility to balancing markets do not perform that well, since at one point in time, they simply have no opportunity to sell energy due to the battery size constraints. Besides that, we see that the battery usages of all strategies are similar to each other.}
    \label{fig: fixed strategies}
\end{figure}

\end{document}